\documentclass[a4paper,11pt]{article}
\usepackage[applemac]{inputenc} 
\usepackage{amsfonts,amsmath,amssymb,amscd}
\usepackage[french]{babel}
\usepackage[T1]{fontenc}

\newcommand{\N}{\mathbb{N}}

\newcommand{\Z}{\mathbb{Z}}
\newcommand{\R}{\mathbb{R}}
\newcommand{\Q}{\mathbb{Q}}

\def\le{\leqslant}
\def\ge{\geqslant}

\newtheorem{defi}{d\'efinition}
\newtheorem{rem}[defi]{- Remarque}
\newtheorem{ex}[defi]{- Exemple}

\newtheorem{prop}{Proposition}
\newtheorem{lem}{Lemme}
\newtheorem{theo}{- Th\'eorème}
\newtheorem{cor}{Corollaire}

\title{L'Algèbre tropicale comme algèbre de la caractéristique  1:\\
 Polynômes rationnels et fonctions polynomiales }
\author{Dominique Castella\\
Laboratoire de Mathématiques\\
Université de la Réunion\\
\small{dominique.castella@univ-reunion.fr}}
 \date{29/08/2008}
 
\begin{document}
 
\maketitle

\begin{center}
 Université de La Réunion\\
 Août  2008
 \end{center}
 
 \begin{abstract}
  We continue, in this second article, the study, initied in [Cas], of the algebraic tools which play a role in tropical algebra. We especially examine here  the polynomial algebras over idempotent semi-fields. This work is motivated by the development of tropical geometry which appears to be the algebraic geometry of  tropical algebra. In fact, the most interesting  object is the image of a polynomial algebra in its semi-field of fractions. We can thus  obtain, over   good  semi-fields, the analog of classical correspondences between polynomials, polynomial functions and varieties of zeros...
 For example, we show that  the algebras of polynomial functions over a tropical curve associated to a polynomial P, is, as in classical algebraic geometry, the quotient of the polynomial algebra by the ideal generated by P.

  \end{abstract}

   \vskip1cm
 
 \section{Introduction}
   Nous avons définis dans ([Cas]) un cadre formel pour l'étude des algèbres de types Max-Plus, Min-Plus et plus généralement des algèbres sur les semi-corps commutatif idempotent qui  sont en fait les ''corps'' de caractéristique 1.  \\ 
   Nous développons maintenant la parentée entre les courbes tropicales et les courbes alg\'ebriques habituelles (il s'agit dans les deux cas de l'ensemble des ''zéros''  d'un polyn\^ome \`a deux variables) et obtenons  de nouveaux outils algébriques pour l'étude de ces courbes.\\
   Ce travail est motivé par l'essort    de la géométrie tropicale qui apparaît comme la géométrie algébrique de l'algèbre tropicale.\\
Les points de non différentiabilité des polynômes à plusieurs variables peuvent en effet être vus comme  une généralisation des zéros des polynômes. \\
Les semi-anneaux de polynômes étant intègres mais non simplifiables, ils admettent bien un semi-corps des fractions mais ne s'y injectent pas... L'objet algébrique adéquat pour l'étude des fonctions polynomiales est alors le polynôme rationnel qui est  l' image d'un polynôme dans le semi-corps des fractions rationnelles, ou de façon équivalente la classe des polynômes ayant même image. \\
Nous obtenons pour ces polynômes un théorème de correspondance avec les fonctions polynomiales, valable  sur les bons semi-corps idempotents, en particulier sur le semi-corps des réels max-plus. Nous obtenons de plus une bonne correspondance entre les variétés de zéros et l'arithmétique de ces polynômes.

Rappelons pour commencer quelques définitions introduites et quelques propriétés montrées dans l'article précédent :\\

\subsection{ Quasi-corps}

  Rappelons qu'un monoïde est un ensemble muni d'une loi interne associative, admettant un \'el\'ement neutre.\\
   
    On dira qu'un \'el\'ement $x$ d'un  monoïde $(G,*)$  est {\it quasi-inversible} s'il existe un 
\'el\'ement $y$ de $G$ tel que $  x*y*x=x$ et $y*x*y=y$.  On dira alors que $x$ et $y$ sont quasi-inverses l'un de l'autre.  
 
 Dans le cas commutatif $y$ est alors unique et est appel\'e le {\it quasi-symétrique} de $x$ (on le notera $x^*$). \\

 Deux \'el\'ements $x$ et $y$ d'un monoïde $G$ sont {\it orthogonaux}  (notation  : $x \perp y$)  si $x*z = x$ et $y*z = y$ implique $z=e$, o\`u $e$ d\'esigne l'\'el\'ement neutre de $G$. On dira que $(x_1,x_2) \in G^2$ est une {\it décomposition orthogonale} de $x \in G$ si $x = x_1*x_2$ et $x_1 \perp x_2$. Dans le cas commutatif on notera $x = x_1 \bigoplus x_2$ pour indiquer que $(x_1,x_2)$ est une décomposition orthogonale de $x$.\\

  Un quasi-anneau est un semi-anneau unitaire dont l'élément unité admet un quasi-symétrique.

 Dans la suite on notera $0$ l'\'el\'ement neutre de l'addition d'un quasi-anneau $A$, et $\epsilon$ le quasi-inverse de $1$. On aura donc, pour tout $x \in A$, $x^*=\epsilon x$.\\
 
 Si $A$ est un quasi-anneau, l'ensemble $A[X_i]_{i\in I}$ des polyn\^omes \`a coefficients dans $A$ est aussi un quasi-anneau.\\
 Un quasi-anneau est simplifiable \`a droite si $\forall (x,y,z) \in A^3$, $x*z = y*z \Longrightarrow x=y$.\\
 
 Un {\it quasi-corps} est un quasi-anneau $(K,+,*)$ tel que $(K^*,*)$ soit un groupe  (o\`u $K^*= K-\{0\})$.\\

  \subsection{ Caract\'eristique d'un semi-anneau}
 
 Rappelons qu'un semi-anneau $A$ se définit comme un anneau, en affaiblissant la condition $(A,+)$ est un groupe commutatif en $(A,+) $ est un monoïde commutatif.\\
 
 Soit $A$ un semi-anneau ; on d\'efinit $H \subset \N$ comme l'ensemble des entiers $k$ tels que $k.1 + 1 = 1$. \\
 
 \begin{prop}
 Il existe un unique $ n \in \N$  tel que $H = n \N$. Cet entier  $n$ est appel\'e caract\'eristique de $A$ (notation $car(A)$).\\
 \end{prop}

 \begin{rem} 
 Les semi-anneaux de caract\'eristique 1 sont les semi-anneaux idempotents  (i.e. tels que $x+x=x$ pour tout $x$).
 \end{rem}
 
  On dira qu'un semi-anneau $A$ est de {\it caract\'eristique pure} p, s'il est de caractéristique $p$ et que,  pour tout $x  \not= 0$, $(k+1)x=x$ implique que $p$ divise $k$.\\
  
  \begin{prop}
  
  a) Tout semi-corps admet une caract\'eristique pure.\\
  b) Les semi-corps de caractéristique non nulle sont des quasi-corps.\\
  c) Les quasi-corps ayant une  caract\'eristique $p$ diff\'erente de $1$ sont des corps. Les quasi-corps de caract\'eristique $1$ sont les semi-corps idempotents.\\
 \end{prop}
    \begin{rem} M\^eme si $K$ est un quasi-corps, l'anneau de polyn\^ome $K[X]$ n'est pas simplifiable :
  si $K$ est de caract\'eristique 1,  $(X+1)(X^2+1) = (X+1)(X^2+X+1)$.\\ En particulier, il ne peut donc pas se plonger dans un quasi-corps des fractions.\\
 \end{rem}
 
 \begin{ex}
Nous utiliserons plus particulièrement dans la suite les deux quasi-corps suivants : \\

Le semi-corps à deux éléments, $F_1 = \{0,1\}$, muni de l'addition telle que $0$ soit élément neutre et $1+1=1$, et de la multiplication habituelle, est un quasi-corps de caractéristique 1,  isomorphe à l'ensemble des parties d'un singleton, muni de la réunion et de l'intersection. Il est facile de vérifier que c'est le seul quasi-corps fini de caractéristique 1...

 Le semi-corps des réels max-plus, $T$, sous la version utilisée dans la plupart des applications, $\R \cup \{-\infty\}$ muni de la loi max comme addition et de la loi $+$ comme multiplication, ou dans sa version, plus pratique pour conserver des notations algébriques générales (et plus facile à suivre par des non-spécialistes) $\R_+$ muni de la loi max comme addition et de la multiplication usuelle...

\end{ex}

 Plus généralement, on a la :\\
 
 \begin{prop}
 Si $A$ est un quasi-anneau de caractéristique 1,  simplifiable, pour tout couple $(x,y) \in A^2$ tel que $xy=yx$, et tout entier $n$, on a :
 $$(x+y)^n = x^n + y^n.$$
 \end{prop}
   
 Un quasi-anneau de caractéristique  1 est ordonné par la relation :  $a \le b$ si $a+b = b$.\\
 Un cas particulier très important est celui des quasi-corps dont l'ordre associé est total. C'est en effet le cas de tous les quasi-corps introduits en algèbre et géométrie tropicale. On parlera alors de quasi-corps totalement ordonnés.\\
Réciproquement, tout groupe totalement ordonné apparaît comme le groupe multiplicatif d'un  quasi-corps, l'addition étant donnée par $a+b = \max(a,b)$. Il suffit en fait que le groupe ait une structure de treillis.\\

   Un {\it module \`a gauche} sur un quasi-anneau $A$ est un triplet $(M,+, .)$ o\`u $(M,+)$ est un quasi-groupe, et $ . $ une loi externe de $A \times M$ dans $M$, v\'erifiant les propri\'et\'es suivantes :\\
   $\forall  a \in A$, $\forall b \in A$, $\forall m \in M$, $\forall n \in M$,  $a.(b.m) = a*b.m$, $(a+b).m = a.m+b.m$, $a.(m+n)=a.m+a.n$ et $1.m = m$.\\
  
  Si $M$ est un module  libre de base $B=(e_i)$, $x$ et $y$ appartenant à $M$ sont orthogonaux si et seulement s'ils ont des supports (relativement à $B$)  disjoints.\\

    \subsection{Points singuliers,  *singuliers et zéros}
    
Si $f$ est un morphisme  d'un quasi-anneau $A$ dans un quasi-groupe $B$, on peut d\'efinir  deux notions duales de  point singulier : 

 On dira que $u$ appartenant \`a $A$ est un point {\it singulier } de $f $, ou que $f$ est singulier en $u$,  s'il existe une d\'ecomposition orthogonale de $f$ dans $Mor(A,B)$, $f = f_1 \bigoplus f_2$, 
telle que $f_1(u) = f_2(u)^*$.  \\

On dira que $u$ appartenant \`a $A$ est un zéro, ou un point {\it *singulier } de $f  $, ou que $f$ est *singulier en $u$,  s'il existe une d\'ecomposition orthogonale de $u$ dans $E$, $u = u_1 \bigoplus u_2$, 
telle que $f(u_1) = f(u_2)^*$.  \\

 On dira de même que   $x = (x_i) \in A^{(I)}$ est un {\it zéro}  d'un polyn\^ome $P \in A[X_i]_{i \in I}$, sur un quasi-anneau commutatif $A$, si $P$ est un zéro (i.e.  un point *singulier) pour le morphisme d'évaluation en $x$,  $P \longmapsto P(x)$ (i.e. si l'on peut écrire $P =P_1\bigoplus P_2$, avec $P_1(x) = P_2(x)^*$).\\

 Cette définition s'étend sans difficulté à un point de $B^{(I)}$, où $B$ est une extension commutative de $A$, ou encore à un point d'une extension (non nécessairement commutative) dans le cas à une variable.\\
 En particulier $0 \in  K^n$ est un zéro  de $P \in K[X_i]$ si et seulement si le terme constant de $P$ est nul.\\
 
Les zéros d'un polynôme à une variable $P \in A[X]$, seront encore appelés {\it racines} de ce polynôme $P$.\\

Si $x \in A$ est un point singulier de l'application polynomiale $P$, $x$ est une racine de $P \in A[X]$, mais la réciproque est fausse,  deux applications polynomiales $P_1$ et $P_2$ correspondant à 2 polynômes orthogonaux, n'étant pas, en général, orthogonales.\\
Cependant, sur le corps des réels (max,+), on peut voir facilement que ces deux notions coïncident : \\
il suffit de voir que l'inf de deux applications définies par des monômes distincts est l'application nulle ;
pour cela, on peut remarquer, si $i \not=j$, que  $a_i x^i \le a_j x^j $ pour tout $x \in \R$ implique $a_i = 0$, en faisant tendre $x$ vers $0$ ou $+ \infty$ suivant les cas...\\ 
   
 \section{ Quotients  et localisations}
  
  A tout morphisme de quasi-anneaux de $A$ dans $B$, est associ\'ee de la mani\`ere habituelle une relation d'\'equivalence compatible avec les lois du quasi-anneau, et donc une structure quotient. Ces relations ne sont par contre pas toutes  obtenues \`a partir d'un id\'eal du quasi-anneau (il faut consid\'erer des sous-quasi-modules convenables de $A^2$). Par contre on peut encore, \`a partir d'un 
 id\'eal d'un quasi-anneau commutatif construire un anneau quotient, et ce cas particulier va nous permettre de g\'en\'eraliser la notion de corps des racines d'un polynôme et d'extension algébrique.\\  Plus g\'en\'eralement, on peut d\'efinir le quotient d'un quasi-module par un sous-quasi-module.\\
 Dans cette section nous supposerons pour alléger le texte que tous les quasi-anneaux considérés sont commutatifs.\\
 
 \subsection{Quotient d'un module par un sous-module}
 Soit $M$ un module \`a droite sur un quasi-anneau $A$ et $N$ un sous-module de $M$.\\
 On notera $[M,N]$ l'id\'eal $\{ a\in A \ / \  M.a \subset N \}$.  Pour $m \in N$, on pose $E_m = \{ r \in M \ / \  m+r \in N \}$.\\
 On d\'efinit une relation $R$ sur $M^2$ par $ \forall (m,n) \in M^2$, $mRn$ si $  (m+N) \cap (n+N )
 \not= \emptyset$ et $\forall a \in  [M,N]$, $E_{ma} = E_{na}$.\\
 $R$ s'appelle la {\it congruence modulo N} (notation $x \equiv y (N)$ pour $x R y$).\\
 
 \begin{prop} $R$ est une relation d'\'equivalence compatible avec les lois de $M$. On notera $M/N$ l'ensemble quotient qui est donc muni d'une structure de module \`a droite.\\
 \end{prop}
 
 D'abord, $R$ est bien une relation d'\'equivalence :\\
 La reflexivit\'e et la sym\'etrie \'etant \'evidentes, il suffit de v\'erifier la transitivit\'e et plus pr\'ecisemment que si $  (m+N) \cap (n+N)
 \not= \emptyset$ et $  (n+N) \cap (p+N) 
 \not= \emptyset$, pour $p \in M$, alors $  (m+N) \cap (p+N) 
 \not= \emptyset$ :\\
 Soient donc $r$, $s$, $t$, $u$ dans $N$ tels que $m+r=n+s$ et $n+t = p+u$ ; on a $m+r+t=n+s+t=p+u+s$ d'o\`u le r\'esultat puisque $M$ est stable par $+$.\\
 
 Maintenant soient $m$, $m'$, $n$ et $n'$ dans $M$ tels que $mRm'$ et $nRn'$. Il est encore clair que 
 $  (m+n+N) \cap (m'+n'+N) 
 \not= \emptyset$, et il suffit donc de montrer que $E_{(m+n)a} =  E_{(m'+n')a}$, pour tout $a \in [M,N]$ :\\
 Si $(m+n)a + u$ appartient \`a $N$, $m'a + (na+u) $, puis $n'a + (m'a+u)$ appartiennent aussi \`a $N$, par hypoth\`ese. \\
 Ceci montre bien la compatibilit\'e avec l'addition. Pour la loi externe, c'est imm\'ediat.\\
 
 \begin{rem}
  Dans le cas des modules la relation $R$ est la relation habituelle, la premi\`ere condition impliquant la seconde.\\
  \end{rem}
  
  \subsection{ Quotient d'un quasi-anneau commutatif par un id\'eal}
  
  \begin{prop}
  Dans le cas d'un quasi-anneau commutatif $A$ , le quotient par un id\'eal $I$ est muni d'une structure de quasi-anneau appel\'e quasi-anneau quotient de $A$ par $I$.\\ Si $I$ est propre ce quotient $A/I$ n'est pas trivial.\\
  \end{prop}
  
  Pour le premier point il suffit de voir que la relation est compatible avec la multiplication, ce qui est encore imm\'ediat.\\
  De plus si $1$ est congru à $0$ modulo $I$, il existe un $i \in I$ tel que $1+i $ appartienne à $I$, et la deuxième condition implique alors $1 \in I$.\\

  \begin{ex}
  Soit $K$ un quasi-corps et $L$ une extension de $K$. On supposera ici $L$ de caract\'eristique 1 et totalement ordonn\'e (par la relation $x \le y$ si $x+y = y$). 
  Si $A = K[X] $,  $I_x$, l'ensemble des polyn\^omes ayant $x \in L$ comme racine, est un id\'eal et le quasi-anneau quotient est  isomorphe au sous-quasi-anneau $K[x]$ de $L$.\\
   De m\^eme si $A =  K[X_1, \cdots, X_n]$ et $x \in K^n$, l'ensemble des polyn\^omes dont $x \in L^n$ est un zéro, est un id\'eal, $I_x$, et si $R_x $ est la relation d\'efinie par cet id\'eal on a $P R_x Q $ si et seulement si $P(x) = Q(x)$.\\
  Supposons que $P$ et $Q$ appartiennent \`a $I_x$ et soit $P=P_1+P_2$, $Q=Q_1+Q_2$ des d\'ecompositions orthogonales de $P$ et $Q$ telles que $P_1(x) = P_2(x)$ et $Q_1(x) = Q_2(x)$.\\
 Dans le cas d'un quasi-corps totalement ordonn\'e, on peut supposer que $P_1$ et $Q_1$ sont des mon\^omes et on peut supposer $P_1 \not= Q_1$, sinon le r\'esultat est imm\'ediat. $P_1+ Q_2$ et $P_2+Q_1$ sont alors orthogonaux et $P+Q$ a bien un zéro en $x$. 
\\Si $P \in I_x$ il est facile de voir que, pour tout $Q$, le polyn\^ome $PQ$ a un zéro en $x$.\\

   \end{ex}
  
  
  On dira qu'un id\'eal $I$ d'un quasi-anneau commutatif $A$ est {\it fermé} si pour tout $a \in A$,  \\ $(a+I)\cap I \not= \emptyset \Longrightarrow  a \in I$.\\ 
  On dira qu'il est {\it dense} si, pour tout $a \in A$, $(a+I) \cap I \neq \emptyset$.\\
 
 Si $I$ est un idéal d'un quasi-anneau $A$, il est facile de vérifier que  $ \overline{I} = \{ x \in A \ / \ (x+I) \cap I \neq \emptyset \}$
est un idéal fermé, la {\it clôture} de $I$.\\
De même il est facile de vérifier que  le  {\it c\oe ur} de $I$, $C(I) = \{ x \in  \overline{I}  / \ \forall \alpha \in  I, \forall r \in A,  x\alpha +r  \in I  \Longrightarrow r \in I \}$  
est un idéal de $A$.\\
Il est clair que $I$ est fermé si et seulement si $\overline{I} = I $.\\

\begin{prop}
Soit $A$ un quasi-anneau commutatif et $I$ un idéal de $A$. La classe de $0$ modulo $I$ est égale au c\oe ur de $I$ et est un idéal fermé de $A$.\\
En particulier si $I$ est fermé, la classe de $0$ modulo $I$ est $I$.\\
\end{prop}

Soit $x \in A$ : on a par définition  $x \equiv 0 (I)$ si $(x+I) \cap I \neq \emptyset$ et si pour tout $\alpha \in I$, $x\alpha + r \in I \Longrightarrow r \in I$, c'est à dire exactement si $x \in C(I)$.\\
Comme la congruence modulo $I$ est compatible avec la structure d'anneau, si $x \equiv 0 (I) $ et $x+y \equiv 0 (I) $ on a nécessairement $y \equiv 0 (I)$ ce qui montre que la classe de $0$ est bien fermée.\\
Si I est fermé, il est clair que $I = C(I)$.\\

\begin{rem}
Si $A$ est un anneau, tout idéal est fermé.\\
Si $A = K[X]$ est un quasi-anneau de polynômes sur un quasi-corps commutatif de caractéristique 1, l'idéal $PK[X]$ des multiples d'un polynôme non nul $P$ est dense si et seulement si $P(0) \not= 0$ ; si $P=X$ il est fermé.\\
\end{rem}

   \subsection{ Id\'eaux premiers, fortement premiers}

 Soit $I$  un idéal d'un quasi-anneau commutatif $A$ ;  $I$ est {\it premier} si le quasi-anneau quotient $A/I$ est intègre ;  il est {\it fortement premier} si le quasi-anneau quotient est simplifiable.\\

 \begin{prop} 
 Tout idéal maximal d'un quasi-anneau commutatif est soit fermé soit dense.\\
 
 \end{prop}
 
 Si $I$ est maximal, soit $\overline{I}=I$ et $I$ est fermé, soit $\overline{I} = A$ et $I$ est dense...\\
 
 \begin{prop}
 Un idéal fermé $K$ d'un quasi-anneau commutatif  $A$ est premier si et seulement si \\ $\forall x \in A - K, \ \forall y \in A, \ xy \in K  \Longrightarrow y \in K$.\\
 Un idéal $I$ est premier si et seulement si son c\oe ur $C(I)$ l'est.\\ 
  \end{prop}
 
    
 
   
   \subsection{Corps de fractions d'un quasi-anneau commutatif intègre}
   
   La localisation des monoïdes s'applique au monoïde $A^*$, si $A$ est un quasi-anneau intègre, et il est facile de vérifier que comme dans la cas classique, on obtient ainsi un quasi-corps  $B = Frac(A)$, que nous appellerons { \it quasi-corps des fractions} de $A$ et un morphisme 
   $i$ (non injectif en général) de $A$ dans $B$, tels que $B$ ne contienne que les produits d'éléments de $i(A)$ et de leurs inverses.\\
   Plus précisément, on a $i(a) = i(b)$ si et seulement si il existe $z \in A^*$ tel que $za = zb$.\\
   
   $i$ est donc injectif si et seulement si $A$ est simplifiable et $i(A)$ est dons toujours un quasi-anneau simplifiable. De plus tout morphisme de $A$ dans un quasi-anneau simplifiable $C$ se factorise en un morphisme de $i(A)$ dans $C$.\\
   On dira donc que $i(A)$ est { \it l'enveloppe simplifiable} de $A$.\\
   
   Ceci s'applique en particulier aux quasi-anneaux de polynômes sur un quasi-anneau intègre $A[X]$.
   Le quasi-corps des fractions $Frac(A[X])$ sera donc isomorphe à l'enveloppe simplifiable de $K[X]$, où $K$ 
   est le quasi-corps des fractions de $A$.\\
   Nous noterons dans la suite $A\{X\}$ l'enveloppe simplifiable du quasi-anneau des polynômes $A[X]$ et $K(X)$ le quasi-corps des fractions de $K[X]$, que nous appellerons  le { \it quasi-corps des fractions rationnelles } à  coefficients dans le quasi-corps $K$.\\
   Dans le cas des anneaux, nous retrouvons bien entendu les notions usuelles et on peut identifier un anneau commutatif intègre avec son enveloppe simplifiable.\\ 
   
   Dans la suite nous appellerons {\it polynômes rationnels} les éléments de $A\{X\}$, puisqu'ils s'identifient aux fractions rationnelles dont le  dénominateur est constant...
   
   \begin{prop}
   Si $P \in K[X_1, \cdots, X_n]$, où $K$est un quasi-corps  commutatif totalement ordonné,  admet un zéro $x = (x_i) \in K^n$ et si, pour $Q \in K[X_1, \cdots, X_n]$, il  existe 
   $R \not= 0$ dans $K[X_1, \cdots, X_n]$ tel que $RP = RQ$,  $x$ est aussi un zéro de $Q$.\\
   \end{prop}

Supposons d'abord $R(x) \not= 0$ et $P(x) \not= 0$ ; il est facile de voir que, s'il y a exactement  k monômes de $R$ prenant la valeur $R(x)$, le nombre de monômes de $RP$ prenant la valeur maximale est au moins $k+1$, ce qui implique que $Q$ est singulier en $x$.\\
Si $R(x) \not= 0$ et $P(x) = 0$, on a bien nécessairement $Q(x) = 0$.\\
   Supposons maintenant $R(x) = 0$. Pour $n = 1$, le résultat est immédiat par simplification.\\
  On procède donc par récurrence sur $n$: \\
  quitte à permuter les variables on peut supposer $x = (x_i)$ avec $x_1 = 0$ ; en écrivant $R = \sum R_i X^i$, $P = \sum P_i X^i$ et $Q = \sum Q_i X^i$, on obtient, si $k$ et $l$ sont les valuations en $i$ de $R$ et de $P$,  $R_k P_l = R_k Q_l$, $P$ étant singulier en $(x_2, \cdots, x_n)$ ; l'hypothèse de récurrence donne bien alors que $Q$ est singulier en $x$.\\
  
      Ceci permet donc de définir la notion de { \it racines d'un polynôme rationnel} $P \in K\{X\}$,  comme étant les racines d'un des représentant et plus généralement, de zéro d'un polynôme rationnel. On peut aussi clairement parler du degré d'un tel polynôme.\\
   Il est clair qu'un polynôme de degré $n$ a au plus $C_{n+1}^2$ racines... En fait il en a au plus $n$ comme nous le verrons ci-dessous.\\
   
  \section{Polynômes et fonctions polynomiales}
  
  $K$ désigne dans la suite un semi-corps idempotent (ou quasi-corps de caractéristique 1) totalement ordonné.\\
  
   \subsection{ Extensions de Quasi-corps}
   
   Le but de cette section est de construire une extension  d'un quasi-corps de caractéristique 1, totalement ordonnée, $K$, contenant une racine d'un polynôme donné de $K[X]$.\\
   
   Pour ceci nous allons montrer que $K\{X\}/(X^n+a)$  est un quasi-corps de caractéristique 1, totalement ordonné, contenant une racine n-ième de a, la classe de $X$.\\
   
  On dira qu'un quasi-corps $K$ est algébriquement clos si tout polynôme de $K[X]$, non constant,  admet au moins une racine dans $K$.\\ 
     
\begin{theo}
Soit $K$ un quasi-corps totalement ordonné de caractéristique 1.\\
a) L'ensemble des polynômes admettant $x$ comme racine est l'idéal $ J = \sum_k (X^{k} + x^k)$.\\
b) L'ensemble des polynômes rationnels de $A =K\{X\}$ admettant $x \in K$ comme racine est  l'idéal  de $A$ engendré par $X+x$.\\
\end{theo}

a) Supposons que $P \in K[X]$ admet $x$ comme racine et
soient $a_iX^i $ et $a_jX^j$ deux monômes distincts, $i <j$,  tels que $P(x) = a_ix^i = a_jx^j$ :\\
on a donc pour tout $k$, $a_kx^k \le P(x)$  et $ P_1 = a_jX^j+a_iX^i = a_jX^i(X^{j-i} +x^{j-i}) $ appartient à $J$ ; 
si $k \ge i$, on pose  $R_k  = a_kx^{k-j}X^j$ et $a_kX^k+a_kx^{k-j}X^j$ appartient aussi à $J$ ; de plus $P = P+R_k$ puisque 
$a_kx^{k-j} \le a_j$ par hypothèse ; si $k \le i$,
$ a_kx^{k-i}X^k(X^{i -k}+ x^{i-k})$ appartient là encore  à $J$ et en posant $R_k  = a_kx^{k-i}X^i$, on a de même $P = P + R_k$.
On obtient ainsi que $ P = P + \sum R_k$ appartient à $J$.\\
La réciproque est  claire.\\

b) Il en découle aussi que  tout multiple de  $X+x$ admet $x$ comme racine ; réciproquement il suffit de montrer que tous les $X^k+x^k$ sont multiple de $X+x$ ; ceci vient de l'égalité 
$$(X+x)^k(X^k+x^k)=\sum_0^{2k} x^iX^{2k-i}= (X+x)^{2k}$$ qui implique dans $K\{X\}$, l'égalité $(X+x)^k=X^k+x^k$.\\

En raisonnant par récurrence sur le degré on a immédiatement le :\\

\begin{cor}
Si $K$ est un quasi-corps totalement ordonné de caractéristique 1, et $P \in K[X]$ un polynôme de degré $n$, $P$ a au plus $ n$ racines dans $K$.
\end{cor}

\begin{theo}
Soit $K$ un quasi-corps totalement ordonné et $a \in K$, $ n\in \N^*$ ;  \\
 $L = K\{X\}/(X^n+a)$ est un quasi corps totalement ordonné contenant $K$ et la classe $x$ de $X$ dans $L$ vérifie $x^n = a $.
\end{theo}

Considérons l'application linéaire $\phi$ de $K[X]$ dans $E = K+KX+KX^2+ \cdots + KX^{n-1}$ définie par $\phi(X^m) = a^qX^r$ si $m=nq+r$ avec $0 \le r < n$. On note $\theta$ le morphisme canonique de $K[X]$ dans $K\{X\}$.\\ Comme le degré est indépendant du représentant, $F = \theta(E)$ est l'ensemble des polynômes rationnels de degré inférieur ou égal à $n-1$. \\ 

\begin{lem}
$P \in K[X]$ est *singulier pour $\phi$ si et seulement si $P \in J = \sum_k (X^{nk} + a^k)$ :\\
\end{lem}

Tout polynôme $P$ peut s'écrire de façon unique $P = \sum_0^{n-1} P_iX^i$ les $P_i$ étant des polynômes en $Y = X^n$, et on a alors
$\phi(P) = \sum_0^{n-1} P_i(a)X^i$.\\
$P$ est *singulier pour $\phi$ si et seulement  chacun des polynôme $P_i$ est  *singulier en $a$, ou encore si $a$ est une racine de $P_i$.\\
Or, d'après la proposition précédente,  $P_i$ admet une racine  en $a$ si et seulement si $P_i \in I = \sum_k (Y^{k} + a^k)$, ce qui donne bien le résultat annoncé.\\

\begin{lem}
Si $P$ et $Q$  ont même image par $\phi$, ils  sont  congrus modulo $J$  :\\
\end{lem}

Il suffit de voir que si $\phi(P) = \phi(Q)$ alors, pour tout $U \in J$ et pour tout $V \in K[X]$, $PU+V$ est 
*singulier pour $\phi$ si et seulement si $QU+V$ l'est ; or $\phi$ est *-singulier en $PU+V$ si et seulement s'il existe une décomposition orthogonale de $V$, $V=V_1+V_2$ telle que $\phi$ soit 
*singulier en $V_2$ et $\phi(V_1) \le \phi(PU)$. \\
L'équivalence annoncée provient alors de ce que $\phi(PU) = \phi(\phi(P)\phi(U)) =\phi(QU)$.\\

On en déduit immédiatement    que dans $A = K[X]/J$ la classe de $X^n$ est égale à la classe de $a$ et 
que $K$ s'injecte dans $A$.  De plus  $A$ est intègre, car le coeur de $J$ est réduit à $\{0\}$ (si $P \not= 0$, il existe un monôme non nul $V$, tel que  $V \le P(X^n+a)$ et $P(X^n+a) + V \in J$, alors que $V$ n'appartient pas à $J$).\\
On vérifie de plus  bien aisément que $K$ s'injecte alors dans le corps des fractions $L$ de $A$ et l'image de $X$, $x$,  y est naturellement  une racine n-ième de a.\\
Pour tout $c \in L$ et tout entier $k$,  on a $(c+x^k)^n = c^n+ \cdots + x^{kn} = c^n+x^{kn}$ et si $c^n \le a^k$, $(x^k)^n \ge c(x^k)^{n-1}$ et donc $x^k \ge c$ ; de même si $a^k \le c^n$, on a $c \ge x^k$.  En particulier $x\ ge 1$ si $a \ge 1$ et $x \le 1$ si $a \le 1$.\\
Il est donc clair que $A$ et  $L$ sont totalement ordonnés ; pour $P \in K[X]$, $P(x)$ est alors égal au monôme dominant évalué en $x$, $a_ix^i$ qui est inversible dans $A$ puisque $x$ l'est. Ceci montre qu'en fait $A$ est un corps et donc égal à $L$.\\  
Soit $\psi$ le morphisme de $K[X]$ dans $L$, $P \longmapsto P(x)$ :  $J$ est l'ensemble des éléments *singuliers pour $\psi$ et $(X^n+a)$ est l'idéal des polynômes rationnels *singuliers  pour le morphisme quotient $\Psi$ de $K\{X\}$ dans $L$. \\ On obtient ainsi un morphisme surjectif de $K[X]$  dans 
$K\{X\}/(X^n+a)$ qui se factorise en un isomorphisme de $A$ dans $K\{X\}/(X^n+a)$ 
On peut ainsi sans inconvénient identifier ces deux quotients.\\

Il est clair que ceci permet de construire un quasi-corps contenant $K$ et toutes les racines d'un polynôme $P$ donné et plus généralement un quasi corps algébriquement clos contenant un quasi-corps totalement ordonné, donné.\\ 
 Sur ce corps $L$, un polynôme $P$ se factorise donc (dans $L\{X\}$) en un produit \\ 
 $\alpha (X+a_1) \cdots (X+a_n)$, où les $a_i$ sont les racines de $P$. \\
  

 
   \subsection{Polynômes rationnels et fonctions polynomiales sur $K^n$}
   
   On définit, comme dans le cas d'une variable, $K\{X_1, \cdots,X_n\}$, l'image de $K[X_1, \cdots, X_n]$ dans son corps des fractions, comme étant le quasi-anneau simplifiable des polynômes rationnels à n variables sur $K$.\\
   Pour un polynôme $P \in K\{X_i\}_{i \in \N}$ on notera  $V(P)$ l'ensemble des zéros de $P$ dans $K^n$. Pour $x \in K^n$ on notera $I_x$ l'ensemble  $\{ Q \in K\{X_i\} \  / \  x $ est un zéro de $Q \}$, des polynômes dont $x$ est un zéro.\\
   
   Soit $P \in K[X_1, \cdots, X_n]$, $ P = \sum_I \lambda_\alpha X^\alpha$, où $I$ est une partie finie de $\N^n$.\\
    On dira que $P$ est {\it convexe} si, pour tout  tout $\gamma \in \N^n$ appartenant à l'enveloppe convexe de $I$, tel que $\gamma = 1/m(\sum_I \gamma_i i) $, on a :  $\lambda_\gamma = \Pi \lambda_i^{\frac{\gamma_i}{m}}$. \\
    Si $K$ est algébriquement clos, on appellera {\it enveloppe convexe } de $P$ le plus petit  polynôme convexe supérieur à $P$.\\   On notera $conv(P)$ cette enveloppe convexe.\\
  On dira qu'un monôme $\lambda_\alpha X^\alpha$ de $P$
   est extrémal  si l'enveloppe convexe du polynôme obtenu en supprimant ce monôme strictement inférieure à celle de $P$. Il est clair que l'enveloppe convexe de $P$ est égale à l'enveloppe convexe de la somme de ses monômes extrémaux.\\
   
   Pour un polynôme $P$ de $K[X_1, \cdots, X_n]$, on notera  $\overline{P}$ le polynôme rationnel   et 
      $\tilde P$ la fonction polynomiale associés.\\
   
   On a alors le lemme (bien connu lorsque $K$ est le quasi-corps des réels max-plus) :
   
   \begin{lem}
   a) La fonction polynomiale définie par $P$ et celle définie par son enveloppe convexe sont les mêmes.\\
  b)  Les fonctions polynomiales $\tilde P$ et $\tilde Q$ sont égales si et seulement si les enveloppes convexes de $P$ et $Q$ sont égales.
   \end{lem}

 a)  Soit $ P = \sum_I \lambda_\alpha X^\alpha$, où $I$ est une partie finie de $\N^n$.\\
   On dira que le monôme  $ \lambda_\alpha X^\alpha$ est localement dominant si \\
   $A_\alpha = \{ x \in K^n \ / \  \lambda_\alpha x^\alpha > \lambda_\beta X^\beta, \forall \beta \in I-\{\alpha\} \ \}$ est non vide.\\
   Il suffit  donc de voir que les monômes localement dominant sont exactement les monômes extrémaux  pour obtenir le résultat.\\
   Il est clair qu'un monôme non extrémal ne peut être dominant. Il reste donc à voir que si $ \lambda_\alpha X^\alpha $ est extrémal il est localement dominant.
   On peut distinguer deux cas : \\
   - soit $ \alpha $ n'est pas dans l'enveloppe  convexe   des $\beta \in I$, 
   $\beta \not= \alpha$.\\
    On peut alors  trouver une $\Q$-forme linéaire sur , $l$, telle que $l(\gamma) \in \Z$ pour tous les $\gamma$ appartenant à $I$ et telle que $l(\alpha) > 0$,  $l(\beta)< 0$ pour tous les $\beta \not= \alpha$ de $I$. Si  $l((\gamma_1, \cdots, \gamma_n)) = \sum c_i\gamma_i$,  il suffit de choisir $x = (x^{c_1}, \cdots, x^{c_n})$, où $x \in K$  est suffisamment grand, ce qui est possible car $K$ est infini et totalement ordonné et ne peut donc avoir de plus grand élément.\\
  - soit $m \alpha = \sum_{I-\{\alpha\}} c_\beta \beta$ avec $m \in \N$, $c_\beta \in \N$,  et $\lambda_\alpha^m > \Pi \lambda_\beta^{c_\beta}$ ;  $ x = (1, \cdots, 1)$ convient alors.\\
  
  b) Si les enveloppes convexes sont différentes, $conv(P+Q) > conv(P)$ ou $conv(P+Q)>conv(Q)$.
  Il existe donc un monôme extrémal de $conv(P+Q)$, n'apparaissant pas, par exemple, dans $conv(P)$.
  Or ce monôme est localement dominant ; soit  $x$ pour lequel ce monôme domine, $(P+Q)(x) $ ne peut donc être égal à $P(x)$ ce qui prouve bien que $\tilde P \not= \tilde Q$.\\
 La réciproque est immédiate.\\

   On peut en déduire qu'il y a une bonne correspondance entre les fonctions polynomiales et les polynômes rationnels :\\

   \begin{theo}
   Si $ K$ est un quasi-corps algébriquement clos infini, l'application naturelle de $K\{X_1, \cdots,X_n\}$ dans les fonctions polynomiales de $K^n$ dans $K$, est un isomorphisme de quasi-anneaux.\\
   \end{theo}
   
Considérons  l'application $\Phi$ de $K[X_1, \cdots, X_n]$ dans l'ensemble des fonctions polynomiales sur $K^n$, $Pol_n(K)$, $P \longmapsto (x \mapsto P(x))$. Si $R$ est un polynôme non nul et $x \in K^n$ est tel que $R(x) \not= 0$,  $R(x)P(x) = R(x)Q(x)$ implique $P(x) = Q(x)$. 
Supposons  $R(x) = 0$, avec $x = (x_i) \in K^n$ ; on peut alors  écrire $ R = \sum_I X_i^jR_i$ où $I$ est  inclus dans l'ensemble des $i \in [1,n]$ tels que $ x_i = 0$ et $j \ge 1$ est tel que $R_i(x) \not=0$ . \\
Choisissons un $i \in I$ : 
en dérivant j fois par rapport à $X_i$, on obtient  $R_i(x)P(x) = R_i(x)Q(x)$ et donc bien $P(x) = Q(x)$.\\
 (La dérivée $P'$ d'un polynôme $\sum a_i Y^i \in K[Y] $ est le polynôme $\sum a_iY^{i-1}  $ et $P \longmapsto P'$ est encore une dérivation).\\
 
   On peut donc bien considérer l'application $\Psi$ de $K\{X_1, \cdots, X_n\}$ dans $Pol_n(K)$ qui associe à un polynôme rationnel la valeur de l'un de ses représentants en $x$ et $\Psi$ est clairement un morphisme surjectif de quasi-anneaux.\\
   Il reste donc à montrer l'injectivité de ce morphisme.\\
   
   Or celle-ci découle de la proposition  suivante :\\
   
   \begin{prop} 
   Soient $P$ et $Q$ deux polynômes sur un quasi-corps $K$ idempotent et totalement ordonné, infini et algébriquement clos.\\
   Les propositions suivantes sont équivalentes :\\
a) $P$ et $Q$ définissent les mêmes polynômes rationnels.\\
b) $P$ et $Q$ ont même enveloppe convexe.\\
c) $P$ et $Q$ définissent les mêmes fonctions polynomiales.\\
\end{prop}

On a vu ci-dessus l'équivalence entre les assertions b) et c) et que a) implique c).\\
Il reste donc à voir que b) impliquent  a) ce qui revient en fait à montrer que $\overline{conv(P)} = \overline{P}$, c'est à dire que les classes de $P$ et de son enveloppe convexe sont les mêmes.\\
Soit donc  $ P = \sum_I \lambda_\alpha X^\alpha$, où $I$ est une partie finie de $\N^n$.\\
Soit $J$ l'enveloppe convexe de $I$ dans $N^n$ et
$\gamma \in J - I$  , tel que $\gamma = 1/m(\sum_I \gamma_i i )$, et tel que  $\lambda_\gamma = \Pi \lambda_i^{\frac{\gamma_i}{m}}$. \\
On obtient, en développant $P^m$,
$(\lambda_\gamma X^\gamma)^m \le P^m$ et ceci montre que :\\
 $\overline{P^m} \ge \overline{conv(P)^m} = \sum_J \lambda_\beta^m  X^{m\beta}$ puisque $K\{X_1, \cdots, X_n\}$ est simplifiable ; ceci donne donc l'égalité car l'autre inégalité est évidente.\\
 Le lemme suivant termine alors la démonstration :\\
 
 \begin{lem}
 Soient $P$ et $Q$ deux polynômes rationnels. S'il existe $m > 0$ tel que $P^m = Q^m$, alors
 $P = Q$.
 \end{lem}
 
 On peut supposer $P$ et $Q$ non nuls.\\
 On a $P^m = P^m+Q^m =(P+Q)^m$ et donc $(P+Q)^{2m} \ge P^{2m} + P^{2m-1}Q \ge P^{2m}$ donne 
 $P^{2m} = P^{2m-1}(P+Q)$ d'où finalement $P = P+Q$, après simplification.\\

\begin{rem}
Un polynôme rationnel $P$ admet donc deux représentants particuliers, l'un maximal, l'enveloppe convexe commune de ses représentants, l'autre minimal la somme des monômes extémaux de l'un de ses représentants.
\end{rem}

 On peut maintenant obtenir, du moins sur le corps $T$ des réels max-plus, le lien habituel entre "variété des zéros" et divisibilité :\\
 
\begin{theo}
Soit $T$ le quasi-corps des réels max-plus.
Si $P$ et $Q$ appartiennent  à $ T\{X_i\}_{i \in \N}$, $V(P) \subset V(Q)$ si et seulement si $P$ divise une puissance de $Q$.\\
\end{theo}

  Pour ceci nous considérerons les parties suivantes de $K^n$ associées à un polynôme rationnel $P =  \sum_I \lambda_\alpha X^\alpha \in T\{X_1, \cdots, X_n\}$ :\\
   $A_\alpha(P) = \{ x \in T^n \ / \  \lambda_\alpha x^\alpha > \lambda_\beta X^\beta, \forall \beta \in I-\{\alpha\} \ \}$.\\
    $B_\alpha(P) = \{ x \in T^n \ / \  \lambda_\alpha x^\alpha \ge \lambda_\beta X^\beta, \forall \beta \in I-\{\alpha\} \ \}$.\\
     $B_{\alpha,\beta} = B_\alpha \cap B_\beta$, pour $(\alpha,\beta) \in I^2$.\\
On peut remarquer que $V(P) $ est la réunion des $B_{\alpha,\beta} $ non vides.\\

\begin{lem}
a) V(P) est la réunion des $B_{\alpha,\beta} $ non vides pour les $\alpha$ et $\beta$ tels que $A_\alpha$ et $A_\beta$ soient eux mêmes non vides.\\
b) Les $A_\alpha$ non vides sont connexes par arcs.\\
c) Le graphe $\Gamma$ dont les sommets sont les $\alpha$ tels que  $A_\alpha$ soit non vide et les arêtes les $B_{\alpha,\beta} $ de codimension 1,  est connexe.\\
\end{lem}
 
 \noindent
 a) Ceci découle de ce qu'au voisinage d'un point de $T^n - V(P)$ est dense dans $K^n$.\\
 b) Ce sont même des convexes sur le corps des réels max-plus... \\
 c)  Si $x \in A_\alpha$ et $y \in A_\beta$, le segment $[x,y]$ a une intersection non vide avec une suite finie de $A_{\alpha_i}$ et coupe nécessairement les $B_{\alpha_i,\alpha_{i+1}}$ ; quitte à prendre une succession de segments pour joindre $x$ à $y$, on peut supposer que ce chemin  ne rencontre que des $B_{\alpha_i,\alpha_{i+1}}$ de codimension 1.  Ceci  fournit alors bien un chemin dans le graphe entre les sommets  $\alpha$ et $ \beta$.\\

\begin{lem}
Si $V(P) = V(Q)$, $\{\alpha \ / \  A_\alpha(P) \not= \emptyset \} = \{ \ \beta \ / \ A_\beta(Q)\not= \emptyset \} $ et les graphes associés sont isomorphes.\\
On peut donc numéroter les monômes  extrémaux de $P$ et de $Q$ de sorte que $A_{\alpha_i}(P) = A_{\beta_i}(Q) = A_i$ pour tout $i$.\\  

\end{lem}

Soient $x, y  \in \{A_\alpha(P)\}$ tels que $x \in A_\beta(Q)$ et $y \in A_\gamma(Q)$ : si $\beta \not= \gamma$, il existe un chemin continu $\phi$ de $x$ à $y$ dans $\{A_\alpha(P)\}$ et ce chemin rencontre nécessaire $V(Q)$, ce qui est contradictoire avec la définition de $\{A_\alpha(P)\}$ et l'hypothèse $V(P) = V(Q)$.\\
On trouve donc pour chaque $\alpha$ un $\beta $ tel que $A_\alpha(P) = A_\beta(Q)$ et par symétrie, en utilisant le a) du lemme précédent, on obtient bien  une bijection entre les sommets du graphe qui respecte les arêtes...

\begin{lem}
Soient les fractions rationnelles $L_{i,j} = \frac{\lambda_{\alpha_i}}{\lambda_{\alpha_j}} X^{\alpha_i- \alpha_j}$ et $M_{i,j} = \frac{\mu_{\beta_i}}{\mu_{\beta_j}} X^{\beta_i - \beta_j}$.\\
a) Si $B_{i,j} $ est de codimension 1, il existe un rationnel  positif $k_{i,j}$   tel que $L_{i,j}= M_{i,j}^{k_{i,j}}$.\\
b) Il existe un entier $m$ tel que, pour tout $x \in A_i$  et tout $j$ :
$L_{i,j}^{m}(x)  \ge M_{i,j}(x)$.\\
\end{lem} 

\noindent
a) Dans le cas des réels max-plus, $B_i$, l'adhérence de $A_i$,  est un convexe défini par le  système d'inéquations  affines (à coefficients entiers), pour tout $r\not=i$, $L_{i,r}(x) \ge 1$, mais aussi par le système d'inéquations, pour tout $r\not=i$, $M_{i,r}(x) \ge 1$. \\ 
Chacune des  équations  $L_{i,j}(x) = 1$ et $M_{i,j}(x) = 1$ définit un hyperplan contenant l'intersection $B_i \cap B_j$ ; mais cet hyperplan  est  unique si $B_{i,j} $ est de dimension $n-1$ et il en découle dans ce cas que  ces deux équations sont proportionnelles, ce qui donne, en notation algébrique, l'existence d'un rationnel $k_{i,j}$ tel que $L_{i,j} = M_{i,j}^{k_{i,j}}$.  Ce rationnel est nécessairement strictement positif car $A_i $ est inclus dans le demi-plan positif pour les deux équations.\\ 
b) Si $A_i \cap A_j$ est vide ou  de codimension strictement plus grande que 1, on peut trouver, d'après le lemme précédent,   dans le graphe $\Gamma $ associé à la fois à $P$ et à $Q$, (en notant $s$ le sommet $\alpha_s = \beta_s$), un chemin allant du sommet $i_0 = i $ au sommet $i_k = j$ ; il existe donc d'après le a),  pour chaque entier $1 \le s \le k$, un rationnel positif $r_s$ tel que $L_{i_{s-1},{i_s}} = M_{i_{s-1},{i_s}}^{r_s}$. \\
 On vérifie, par récurrence sur la distance $k$, de $i$ à $j$ dans le graphe $\Gamma$, qu'il existe des rationnels positifs $r$ et $s$ tels que, pour tout $x \in A_i$ : \\
 $M_{i,j}^r(x) \le L_{i,j}(x)  \le  M_{i,j}^t(x) $.\\
 Pour $k = 1$ c'est une conséquence immédiate du a).\\
 Si l'on suppose  le résultat vrai au rang $k$ et que la distance de $i$à $j $ est égale à $k+1$,  il existe un sommet $l$ tel que $B_{l,j}$ soit de codimension 1, et des rationnels positifs $r$ et $t$ tels que, pour tout $x \in A_i$,  $M_{i,l}^r(x)  \le L_{i,l}(x) \le  M_{i,l}^t(x)  $.\\
 On a  $L_{i,j} = L_{i,l}L_{l,j}$ et $M_{i,j} = M_{i,l}M_{l,j}$ ; pour $x \in A_i$, $L_{i,l}(x) \ge 1$ et $M_{i,l}(x) \ge 1$.\\
  Si $L_{l,j}(x)$ est plus grand que 1, et on obtient immédiatement 
 $\max(r,t) M_{i,j}(x)  \ge L_{i,j} \ge \inf(r,t) M_{i,j}(x)$.\\ 
 Si $L_{l,j}(x)$ est plus petit que 1, il en est de même de $M_{l,j}(x)$ d'après le a) ; on a alors $r M_{i,j}(x) \ge L_{i,j}(x)$, et aussi $\frac{1}{r} M_{i,j}(x) \ge L_{i,j}(x)$.\\ 
 
 La récurrence est donc vérifiée et il existe donc des rationnels $t_{i,j}$ tels que, pour tout $x \in A_i$, 
 $L_{i,j}(x) \ge M_{i,j}^{t_{i,j}}(x)$.\\
 On peut alors  choisir un entier non nul $m_j$ tel que $m_j t_{i,j}$ soit entier et on a alors pour tout $x \in A_i$, $L_{i,j}^{m_i}(x)  \ge M_{i,j}(x)$ puisque $M_{i,j}(x)  \ge 1$.\\ 

 Les $L_{i,j}(x)$ étant supérieurs à 1 sur $A_i$, le plus grand des $m_i$ convient.\\

\begin{lem}
Si $V(P) = V(Q)$, il existe un entier non nul $k$  tels que $Q$ divise  $P^{k}$.
 \end{lem} 

Pour un $\gamma \ge \sup_i(\beta_i)$,
considérons le polynôme $R_k = \sum_r \frac{\lambda_{\alpha_r}^k}{\mu_{\beta_r}} X^{k\alpha_r - \beta_r + \gamma}$.\\
Grâce au lemme précédent, on peut choisir $k$ assez grand pour que, pour chaque i,   le monôme $\frac{\lambda_{\alpha_i}^k}{\mu_{\beta_i}} X^{k\alpha_i - \beta_i+\gamma}$ soit dominant sur $A_i$.\\
On obtient alors pour tout $x \in \cup A_i$ :\\
$P^k(x)x^\gamma = \lambda_{\alpha_i} x^{\alpha_i+\gamma}=R(x)Q(x)$. Cette égalité est alors vraie par densité pour tout $x\in T^n$.
Le théorème précédent donne donc $P^kX^\gamma = QR$.\\
Mais pour les $i$ tels que $\gamma_i \not= 0$, soit  $X_i$ ne divise pas $Q$ et nécessairement $X_i^{\gamma_i} $ divise $R$, soit $X_i$ divise $Q$ et la condition $V(Q) = V(P) $ prouve que $X_i$ divise $P$. Quitte à changer $k$ on peut supposer dans tous les cas  que $Q$ divise $P^k$.\\

Soient maintenant $P$ et $Q$ tels que $V(P) \subset V(Q)$. On a $V(PQ) = V(Q)$ et il existe donc $k$ tel que $(PQ)$ divise $ Q^{k}$, ce qui donne le résultat annoncé.\\

On obtient alors facilement la :\\

\begin{prop}
Soit $T$ le quasi-corps des réels max-plus et 
soit $I$ l'idéal engendré par un polynôme $P \in T\{X_1,\cdots, X_n\}$. \\
a) Les classes de deux polynômes $A$ et $B$ de  $T\{X_1,\cdots, X_n\}$. modulo $I$, sont égales si et seulement si $A(x) = B(x)$ pour tous les  $x \in V(P)$ tels que $P(x) \not=0$.\\
b) Le radical de $I$,  $rad I = \{Q \in  T \{ X_1,\cdots, X_n\} \ / \ \exists k \in \N, \  Q^k \in I \}$,  est égal  à l'intersection de tous les    $Ker \  \epsilon_x$ pour $x \in V(P)$, où   $\epsilon_x $ est le morphisme d'évaluation en $x$, $P \longmapsto P(x)$.\\
\end{prop}

\noindent
a) En effet, $A$ et $B$ sont congrus modulo $I$ si et seulement si, pour tout couple $(U,R) \in I \times  T\{X_1,\cdots, X_n\}$, $AU=R \in I \Leftrightarrow BU+R \in I$. Mais $APU + R$ appartient à $I$ si et seulement si $V(P) \subset V(APU+R) $, c'est à dire si et seulement si, pour chaque  $x \in V(P)$, $R$ est singulier en $x$ ou $R(x) \le APU(x)$.\\
Il est donc facile de voir, que si $A(x) = B(x)$ pour tout $x \in V(P)$, tel que $P(x) \not= 0$, $APU+R \in I $ si et seulement si $BPU+R \in I$.\\
Réciproquement, en considérant les $R \in T$, on voit que la condition est bien nécessaire.\\
b) Ceci découle directement du théorème précédent. \\

    \begin{center}
{\large  REFERENCES}
\end{center}

 [ABG] M. AKIAN, R. BAPAT, S. GAUBERT.{\it Max-plus algebras},  Handbook of Linear Algebra (Discrete Mathematics and Its Applications, L.HOGBEN ed.), Chapter  25, vol. 39, Chapman  \& Hall / CRC, 2006. \\ 
 
 [Cas] D. CASTELLA. L'algèbre tropicale comme algèbre de la caractéristique 1 : Algèbre linéaire sur les semi-corps idempotents. Preprint\\
 
  [IMS] I. ITENBERG, G. MIKHALKIN, E. SHUSTIN. Tropical algebraic geometry. Oberwolfach Seminars, 35, Birkhauser (2007).\\

[ [Izh1] Z. IZHAKIAN. Tropical arithmetic and algebra of tropical matrices. ArXiv:math. AG/0505458, 2005.\\

[Izh2] Z. IZHAKIAN. Tropical varieties, ideals and an algebraic nullstensatz. ArXiv:math. AC/0511059, 2005\\

 [Mik2]  G. MIKHALKIN. Enumerative tropical algebraic geometry in $\R^2$.  J. Amer. Math. Soc. 18 (2005) 313-377.\\
 
[RST] J. RICHTER-GEBERT, B. STURMFELS, T. THEOBLAND. First steps in tropical geometry. Contemporary Mathematics, 377, Amer. Math. Soc. (2005), 289-317.\\

[SI] E. SHUSTIN, Z. IZHAKIAN A tropical nullstensatz. Proc. Amer. Math. Soc. 135 (12), 3815,3821 (2007).\\

\end{document}